\documentclass{amsart}[12pt]
\usepackage[dvips]{graphicx}
\usepackage{amsfonts, amsmath, amssymb, amscd, latexsym, pb-diagram}

\def \Z{{\mathbb Z}}

\newtheorem{thm}{Theorem}[section]    % Standard theorem environment
\newtheorem{corol}[thm]{Corollary}%

\title{Obstructions for subgroups of Thompson's group $V$}

\author{Jos\'e Burillo}
\address{Departament de Matem\`atica Aplicada IV,
Universitat Polit\`ecnica de Catalunya,
Escola Polit\`ecnica Superior de Castelldefels,
Esteve Terrades 5,
08860 Castelldefels (Barcelona), Spain}
\email{burillo@ma4.upc.edu}

\author{Sean Cleary}
\address{Department of Mathematics,
The City College of New York, NY, USA}
\email{cleary@sci.ccny.cuny.edu}

\author{Claas E. R\"over}
\address{School of Mathematics, Statistics and Applied Mathematics,
National University of Ireland Galway,
University Road, Galway, Ireland}
\email{claas.roever@nuigalway.ie}

\begin{document}

\thanks{The authors are grateful for the hospitality of Durham University during the Symposium on Cohomological and Geometric Group Theory.
The first author acknowledges support from MEC grant
MTM2011--25955.  The second author acknowledges support from
the National Science Foundation and that this work was partially
supported by a grant from the Simons Foundation (\#234548 to Sean
Cleary).}

\begin{abstract}Thompson's group $V$ has a rich variety of subgroups, containing all
finite groups, all finitely generated free groups and all finitely
generated abelian groups, the finitary permutation group of a
countable set, as well as many wreath products and other families of
groups.  Here, we describe some obstructions for a given group to be a
subgroup of $V$.
\end{abstract}

\maketitle

Thompson constructed a finitely presented group now known as $V$ as an
early example of a finitely presented infinite simple group.  The
group $V$ contains a remarkable variety of subgroups, such as the
finitary infinite permutation group $S_{\infty}$, and hence all
(countable locally) finite groups, finitely generated free groups,
finitely generated abelian groups, Houghton's groups, copies of
Thompson's groups $F$, $T$ and $V$, and many of their generalizations,
such as the groups $G_{n,r}$ constructed by Higman
\cite{MR0376874}. Moreover, the class of subgroups of $V$ is closed
under direct products and restricted wreath products with finite or
infinite cyclic top group.

In this short survey, we summarize the development of properties of
$V$ focusing on those which prohibit various groups from occurring as
subgroups of $V$.

Thompson's group $V$ has many descriptions.  Here, we simply recall
that $V$ is the group of right-continuous bijections from the unit
interval $[0,1]$ to itself, which map dyadic rational numbers to
dyadic rational numbers, which are differentiable except at finitely
many dyadic rational numbers, and with slopes, when defined, integer
powers of 2. The elements of this group can be described by reduced
tree pair diagrams of the type $(S,T,\pi)$ where $\pi$ is a bijection
between the leaves of the two finite rooted binary trees $S$ and $T$.

Higman \cite{MR0376874} gave a different description of $V$, which he
denoted as $G_{2,1}$ in a family of groups generalizing $V$.

\section{Obstructions}

Higman \cite{MR0376874} described several important properties of $V$
which can serve as obstructions to subgroups occurring in $V$.

\begin{thm}[\cite{MR0376874}]
An element of infinite order in $V$ has only finitely many roots.
\end{thm}

This prevents all Baumslag-Solitar groups $B_{m,n} = \langle a,b\mid
a^nb = ba^m\rangle$ from occurring as subgroups of $V$, if $m$
properly divides $n$; see \cite{CRthesis}.

\begin{thm} [\cite{MR0376874}]
Torsion free abelian subgroups of $V$ are free abelian, and their
centralizers have finite index in their normalizers in $V.$
\end{thm}

This prevents $ GL_n(\mathbb Z)$ from occurring as a subgroup of $V$ for $n\ge 2$.

A group is {\em torsion locally finite} if every torsion subgroup is
locally finite.  That is, if every finitely generated torsion subgroup
is finite.  R\"over \cite{MR1714140} showed

\begin{thm} [\cite{MR1714140}]
Thompson's group $V$ is torsion locally finite.
\end{thm}

This rules out many branch groups from occurring as subgroups of $V$,
including the Grigorchuk groups  of intermediate growth \cite{MR565099} and the
Gupta-Sidki groups \cite{MR696534}. It also rules out Burnside groups.

Holt and R\"over \cite{MR2274726} showed that $V$ has indexed co-word
problem.

\begin{thm}  [\cite{MR2274726}]
The set of words (over an arbitrary but fixed finite generating set)
which do not represent the identity in $V$ is an indexed language, and
hence can be recognized by a nested-stack automaton.
\end{thm}

This property is not easy to verify, however.  But it is inherited by
finitely generated subgroups (see \cite{MR2274726}), and hence groups
which do not have an indexed co-word problem cannot occur as a
subgroup of $V$.

Lehnert and Schweitzer \cite{MR2323454} improved this result.

\begin{thm} [\cite{MR2323454}]
The set of words (over an arbitrary but fixed finite generating set)
which do not represent the identity in $V$ is a context-free language,
and hence can be recognized by a pushdown automaton.
\end{thm}

Again, this property is inherited by finitely generated subgroups, but
the condition is still not easy to verify.

More recently, Bleak and Salaza-D\'iaz \cite{MR3091272} and
subsequently Corwin \cite{Corwin}, using similar techniques showed

\begin{thm} [\cite{MR3091272,Corwin}]
Neither the free product $\Z*\Z^2$ nor the standard restricted
wreath product $\Z\wr\Z^2$ with $\Z^2$ as top group are subgroups of $V$.
\end{thm}

One theorem of Higman \cite{MR0376874} together with a metric estimate
of Birget \cite{MR2104771} gives another obstruction.

\begin{thm} [\cite{MR0376874}]
For any element $v$ of infinite order in $V$, there is a power $v^n$
such that for the reduced tree pair diagram $(S,T,\pi)$ for $v^n$, there
is a leaf $i$ in the source tree $S$ which is paired with a leaf $j$
in the target tree $T$ so that $j$ is a child of of $i$.
\end{thm}

\begin{thm} [\cite{MR2104771}]
For any finite generating set of $V$, There are constants $C$ and $C'$
such that word length $|v|$ of an element of $V$ with respect to that
generating set satisfies $ C n \leq |v| \leq C' n \log{n} $ where $n$
is the size of the reduced tree pair diagram representing $v$.
\end{thm}

Since the powers of $v^n$ will have length thus growing linearly,
these two theorems give as a consequence the following.

\begin{thm}
Cyclic subgroups of $V$ are undistorted.
\end{thm}

We note that this argument applies as well to generalizations of $V$
where there is a linear lower bound on word length in terms of the
number of carets, such as braided versions of $V$ \cite{MR2514382}.

This last theorem has an obvious corollary.

\begin{corol}
If a group embeds in $V$, its cyclic subgroups must be undistorted.
\end{corol}

The reason for this is that in a chain of subgroups $G\supset H\supset
K$ the distortion of $K$ in $H$ cannot be larger than the distortion
of $K$ in $G$.

%For instance, this shows that the presence of distorted cyclic
%subgroups in multiple-slope versions of $F$ such as $F(2,3)$
%\cite{MR2794610} (and clearly thus also multiple-slope versions of
%$V$) means they cannot occur as subgroups of $V$.

We note that Bleak, Bowman, Gordon, Graham, Hughes, Matucci and
J. Sapir \cite{MR3134027} used Brin's methods of revealing pairs for
elements of $V$ to show that cyclic subgroups of $V$ are undistorted.

This result excludes all Baumslag-Solitar groups with $\lvert n\rvert
\neq \lvert m\rvert$, as these have distorted cyclic subgroups. It
also rules out nilpotent groups which are not virtually abelian. An
alternative argument excluding the Baumslag-Solitar groups is due to
Bleak, Matucci and Neunh\"offer \cite{2013arXiv1312.1855B}.

\bibliography{vdistort}

\def\cprime{$'$}
\begin{thebibliography}{10}

\bibitem{MR2104771}
Jean-Camille Birget.
\newblock The groups of {R}ichard {T}hompson and complexity.
\newblock {\em Internat. J. Algebra Comput.}, 14(5-6):569--626, 2004.
\newblock International Conference on Semigroups and Groups in honor of the
  65th birthday of Prof. John Rhodes.

\bibitem{2013arXiv1312.1855B}
C.~{Bleak}, F.~{Matucci}, and M.~{Neunh{\"o}ffer}.
\newblock {Embeddings into Thompson's group $V$ and $co\mathcal{CF}$ groups}.
\newblock {\em arXiv e-prints}, December 2013.

\bibitem{MR3134027}
Collin Bleak, Hannah Bowman, Alison Gordon~Lynch, Garrett Graham, Jacob Hughes,
  Francesco Matucci, and Eugenia Sapir.
\newblock Centralizers in the {R}. {T}hompson group {$V_n$}.
\newblock {\em Groups Geom. Dyn.}, 7(4):821--865, 2013.

\bibitem{MR3091272}
Collin Bleak and Olga Salazar-D{\'{\i}}az.
\newblock Free products in {R}. {T}hompson's group {$V$}.
\newblock {\em Trans. Amer. Math. Soc.}, 365(11):5967--5997, 2013.

\bibitem{MR2514382}
Jos{\'e} Burillo and Sean Cleary.
\newblock Metric properties of braided {T}hompson's groups.
\newblock {\em Indiana Univ. Math. J.}, 58(2):605--615, 2009.

\bibitem{Corwin}
Nathan Corwin.
\newblock {\em Embedding and non embedding results for R. Thompson’s group
  $V$ and related groups}.
\newblock PhD thesis, University of Nebraska \--- Lincoln, 2013.

\bibitem{MR565099}
R.~I. Grigor{\v{c}}uk.
\newblock On {B}urnside's problem on periodic groups.
\newblock {\em Funktsional. Anal. i Prilozhen.}, 14(1):53--54, 1980.

\bibitem{MR696534}
Narain Gupta and Sa{\"{\i}}d Sidki.
\newblock On the {B}urnside problem for periodic groups.
\newblock {\em Math. Z.}, 182(3):385--388, 1983.

\bibitem{MR0376874}
Graham Higman.
\newblock {\em Finitely presented infinite simple groups}.
\newblock Department of Pure Mathematics, Department of Mathematics, I.A.S.
  Australian National University, Canberra, 1974.
\newblock Notes on Pure Mathematics, No. 8 (1974).

\bibitem{MR2274726}
Derek~F. Holt and Claas~E. R{\"o}ver.
\newblock Groups with indexed co-word problem.
\newblock {\em Internat. J. Algebra Comput.}, 16(5):985--1014, 2006.

\bibitem{MR2323454}
J.~Lehnert and P.~Schweitzer.
\newblock The co-word problem for the {H}igman-{T}hompson group is
  context-free.
\newblock {\em Bull. Lond. Math. Soc.}, 39(2):235--241, 2007.

\bibitem{MR1714140}
Claas~E. R{\"o}ver.
\newblock Constructing finitely presented simple groups that contain
  {G}rigorchuk groups.
\newblock {\em J. Algebra}, 220(1):284--313, 1999.

\bibitem{CRthesis}
Claas~E. R{\"o}ver.
\newblock {\em Subgroups of finitely presented simple groups}.
\newblock PhD thesis, University of Oxford, 1999.

\end{thebibliography}

\bibliographystyle{plain}

\end{document}